\documentclass[12pt,a4j,amsfonts]{article}
\usepackage{enumerate,verbatim,calc}
\usepackage{amsmath}
\usepackage{amssymb}
\usepackage{amscd}
\usepackage{ascmac}

\newtheorem{Theorem}{Theorem}
\newtheorem{Lemma}[Theorem]{Lemma}
\newtheorem{Corollary}[Theorem]{Corollary}

\newtheorem{Definition}[Theorem]{Definition}
\newtheorem{Remark}[Theorem]{Remark}
\newtheorem{Example}[Theorem]{Example}

\newtheorem{Fact}[Theorem]{Fact}
\newtheorem{Question}[Theorem]{Question}
\newtheorem{Problem}[Theorem]{Problem}

\newcommand{\R}{\mathbb{R}}

\newcommand{\specialarrow}[3]{\mathrel{\mathop{#1}\limits^{#2}_{#3}}}

\newcommand{\qed}{{\unskip\nobreak\hfil\penalty50\quad\null\nobreak\hfil{\bf
q.e.d.}\parfillskip0pt\finalhyphendemerits0\par\medskip}}

 \DeclareMathOperator{\Spec}{Spec}

 \DeclareMathOperator{\rank}{rank}

\usepackage[dvipdfmx]{graphicx}

\begin{document}

\title{On the limit of Frobenius in the Grothendieck group}
\author{Kazuhiko Kurano (Meiji University) 
\footnote{The author is partially supported by JSPS KAKENHI Grant 24540054. }
\footnote{E-mail: {\sffamily kurano@isc.meiji.ac.jp}}
   \and Kosuke Ohta (Meiji University)
\footnote{E-mail: {\sffamily k\_ohta@meiji.ac.jp}}}
\date{\small Dedicated to Professor Ng\^{o} Vi\d{\^{e}}t Trung for his 60th birthday.}
\maketitle

\begin{abstract}
Considering the Grothendieck group modulo numerical equivalence,
we obtain the finitely generated lattice $\overline{G_0(R)}$
for a Noetherian local ring $R$.
Let $C_{CM}(R)$ be the cone in $\overline{G_0(R)}_{\Bbb R}$
spanned by cycles of maximal Cohen-Macaulay $R$-modules.
We shall define the fundamental class $\overline{\mu_R}$ of $R$
in $\overline{G_0(R)}_{\Bbb R}$, which is the limit of the Frobenius
direct images (divided by their rank) $[{}^e R]/p^{de}$ in the case
${ch}(R) = p > 0$.
The homological conjectures are deeply related to the problems
whether  $\overline{\mu_R}$ is in the Cohen-Macaulay cone  $C_{CM}(R)$ or 
the strictly nef cone $SN(R)$ defined below.
In this paper, we shall prove that  $\overline{\mu_R}$ is in $C_{CM}(R)$
in the case where $R$ is FFRT or F-rational.
\end{abstract}

\section{Introduction}

We shall define the Cohen-Macaulay cone $C_{CM}(R)$,
the strictly nef cone $SN(R)$,
and the fundamental class $\overline{\mu_R}$ for a Noetherian local domain $R$.
They satisfy
\[
\begin{array}{ccccc}
\overline{G_0(R)}_{\Bbb R} & \supset & SN(R) & \supset & C_{CM}(R) - \{ 0 \} \\
\cup & & & & \\
\overline{G_0(R)}_{\Bbb Q} & \ni & \overline{\mu_R}
\end{array}
\]
where ${G_0(R)}$ is the Grothendick group of finitely generated $R$-modules,
$\overline{G_0(R)}$ is the Grothendick group modulo numerical equivalence, and
$\overline{G_0(R)}_{K} = \overline{G_0(R)}\otimes_{\Bbb Z}K$.
By \cite{K23}, $\overline{G_0(R)}$ is a finitely generated free ${\Bbb Z}$-module.
We define $C_{CM}(R)$ to be the cone in $\overline{G_0(R)}_{\Bbb R}$ spanned
by cycles corresponding to maximal Cohen-Macaulay $R$-modules.
If $R$ is F-finite with residue class field algebraically closed,
the fundamental class $\overline{\mu_R}$ is the limit of the Frobenius
direct images (divided by their rank) $[{}^eR]/p^{de}$ as in Remark~\ref{Remark8} (3).
In the case where $R$ contains a regular local ring $S$ such that
$R$ is contained in a Galois extension $B$ of $S$,
then  $\overline{\mu_R}$ is described in terms of $B$ as in 
Remark~\ref{Remark8} (2).

The fundamental class is deeply related to the homological conjectures
as in Fact~\ref{fact8}.
The fundamental class $\overline{\mu_R}$ is in $C_{CM}(R)$ for any
complete local domain $R$ if and only if the small Mac conjecture is true.
Roberts proved $\overline{\mu_R} \in SN(R)$ for any
Noetherian local ring $R$ of characteristic $p > 0$
in order to show the new intersection theorem
in the mixed characteristic case~\cite{Ri}.
In order to extend these results, we are mainly interested in the problem whether $\overline{\mu_R}$ is in such  cones or not.
\begin{Problem}\label{yosou2}
\begin{rm}
If $R$ is an excellent Noetherian local domain, is
$\overline{\mu_R}$ in $C_{CM}(R)$?
\end{rm}
\end{Problem}

Problem \ref{yosou2} is affirmative if  $R$ is a complete intersection. 
However,  even if $R$ is a Gorenstein ring which contains a field, Problem \ref{yosou2} is an open question. 

The following theorem is the main result in this paper. 
We define the terminologies later.

\begin{Theorem}\label{th11}
Assume that $R$ is an F-finite Cohen-Macaulay local domain 
of characteristic $p > 0$
with residue class field algebraically closed.
\begin{itemize}
\item[(1)]
If $R$ is FFRT, then there exist a natural number $n$ and a 
maximal Cohen-Macaulay $R$-module $N$ such that
$n \mu_R = [N]$ in $G_0(R)_{\Bbb Q}$.
In particular,  $\overline{\mu_R}$ is contained in $C_{CM}(R)$.
\item[(2)]
If $R$ is F-rational, then $\overline{\mu_R}$ is contained in $Int(C_{CM}(R))$.
\end{itemize}
\end{Theorem}

In the case FFRT, we shall show that
the cone generated by $[M_1]$, \ldots, $[M_s]$ (in Definition~\ref{FFRT}) 
contains ${\mu_R}$.
In the case of F-rational, the key point  in our proof is to use the
dual F-signature defined by Sannai~\cite{Sa}.

Finally we shall give a corollary (Corollary~\ref{cor14}),
which was first proved in \cite{CK}.

\section{Cohen-Macaulay cone}
In this paper, let $R$ be a $d$-dimensional Noetherian local domain
such that one of the following conditions 
are satisfied:
\begin{itemize}
\item[(a)] 
$R$ is a homomorphic image of an excellent regular local ring containing ${\Bbb Q}$.
\item[(b)]
$R$ is essentially of finite type over a field, ${\Bbb Z}$ or a  complete 
DVR.
\end{itemize}
If either (a) or (b) is satisfied,
there exists a regular alteration of $\Spec R$ by 
de~Jong's theorem~\cite{dJ}.

We always assume that modules are finitely generated.

Let $G_0(R)$ be the Grothendieck group of finitely generated $R$-modules, that is,
\[
G_{0}(R) =
\frac{ \displaystyle\bigoplus_{\mbox{\scriptsize $M$ : f.\ g.\ $R$-module}} {\Bbb Z}[M] }{
<
[M] - [L] - [N] \mid \mbox{$0\rightarrow L \rightarrow M \rightarrow N \rightarrow0$ is exact}
>
} ,
\]
where $[M]$ denotes the generator corresponding to an $R$-module $M$.
Let $C(R)$ be the category of bounded complexes of finitely generated 
$R$-free modules such that every homology is of finite length.
Let $C_d(R)$ be the subcategory of $C(R)$ consisting of complexes of length $d$
with $H_0({\Bbb F}.) \neq 0$.
A complex ${\Bbb F}.$ in $C_d(R)$ is of the form
\[
0 \rightarrow F_{d} \rightarrow F_{d-1} \rightarrow \cdots
\rightarrow F_{1}\rightarrow F_{0}\rightarrow 0 .
\]
For example, the Koszul complex of a parameter ideal belongs to  $C_d(R)$.

For ${\Bbb F}. \in C(R)$, we have a well-defined map
\[
\chi_{{\Bbb F}.} : G_0(R) \longrightarrow {\Bbb Z}
\]
by $\chi_{{\Bbb F}.}([M]) = \sum_i(-1)^i\ell_R(H_i({\Bbb F}.\otimes_RM))$.
We have the induced maps $\chi_{{\Bbb F}.} : G_0(R)_{\Bbb Q} \longrightarrow {\Bbb Q}$
and $\chi_{{\Bbb F}.} : G_0(R)_{\Bbb R} \longrightarrow {\Bbb R}$.
We say that $\alpha \in G_0(R)$ ($\alpha \in G_0(R)_{\Bbb Q}$ or $\alpha \in G_0(R)_{\Bbb R}$) is numerically equivalent to $0$ 
if $\chi_{{\Bbb F}.}(\alpha) = 0$ for any ${\Bbb F}. \in C(R)$.
We define the Grothendieck group modulo numerical equivalence
as follows:
\[
\overline{G_0(R)} = G_0(R)/
\{ \alpha \in G_0(R) \mid 
\mbox{$\chi_{{\Bbb F}.}(\alpha) = 0$ for any ${\Bbb F}. \in C(R)$}
\} .
\]

Then, by Theorem~3.1 and Remark~3.5 in \cite{K23}, 
we know that $\overline{G_0(R)}$ is a non-zero finitely generated
${\Bbb Z}$-free module.\footnote{We need the existence of a regular alteration
in the proof of this result.}

\begin{Example}\label{rei}
\begin{rm}
\begin{itemize}
\item[(1)]
If $d \le 2$, then $\overline{G_0(R)} = {\Bbb Z}$ (Proposition~3.7 in \cite{K23}).
If $d \ge 3$, there exists an example of $d$-dimensional Noetherian local domain $R$ such that $\rank \overline{G_0(R)} = m$ for any positive integer $m$ as in 
(2) (b) (i) below.

\item[(2)]
Let $X$ be a smooth projective variety with embedding $X \hookrightarrow {\Bbb P}^n$.
Let $R$ (resp.\ $D$) be the affine cone (resp.\ the very ample divisor) of this embedding. 
Let $A_*(R)$ be the Chow group of $R$. 
By \cite{K23}, we can define numerical equivalence also on $A_*(R)$, that is compatible with the Riemann-Roch theory as below. 
Let $CH^\cdot(X)$ (resp. $CH^\cdot_{num}(X)$) be the Chow ring (resp. Chow ring modulo numerical equivalence) of $X$. 
It is well-known that $CH^\cdot_{num}(X)_{\Bbb Q}$ is a finite dimensional 
${\Bbb Q}$-vector space.
Then, we have the following commutative diagram:
\[
\begin{array}{ccccc}
G_0(R)_{\Bbb Q} &\specialarrow{\longrightarrow}{\tau_R}{\sim} & A_*(R)_{\Bbb Q} & 
\stackrel{\sim}{\longleftarrow} & CH^\cdot(X)_{\Bbb Q}/ D \cdot CH^\cdot(X)_{\Bbb Q} \\
\downarrow & & \downarrow & & \downarrow \\
\overline{G_0(R)}_{\Bbb Q} & \specialarrow{\longrightarrow}{\overline{\tau_R}}{\sim} & \overline{A_*(R)}_{\Bbb Q} & 
\stackrel{\phi}{\longleftarrow} & CH^\cdot_{num}(X)_{\Bbb Q}/ D \cdot CH^\cdot_{num}(X)_{\Bbb Q}
\end{array}
\]
\begin{itemize}
\item[(a)]
By the commutativity of this diagram, $\phi$ is a surjection.
Therefore, we have
\begin{equation}\label{rank1}
\rank \overline{G_0(R)} \le \dim_{\Bbb Q} CH^\cdot_{num}(X)_{\Bbb Q}/ D \cdot CH^\cdot_{num}(X)_{\Bbb Q} .
\end{equation}
\item[(b)]
If $CH^\cdot(X)_{\Bbb Q} \simeq CH^\cdot_{num}(X)_{\Bbb Q}$, 
then we can prove that $\phi$ is an isomorphism (\cite{K23}, \cite{RS}).
In this case, the equality holds in (\ref{rank1}).
Using it, we can show the following:
\begin{itemize}
\item[(i)]
If $X$ is a blow-up at $n$ points of ${\Bbb P}^k$ ($k \ge 2$),
then $\rank \overline{G_{0}(R)} = n+1$.
\item[(ii)]
If $X = {\Bbb P}^{m} \times {\Bbb P}^{n}$, then 
$\rank \overline{G_{0}(R)} = min\{ m, n \}$.
\end{itemize}
\item[(c)]
There exists an example such that $\phi$ is not an isomorphism \cite{RS}.

Further, Roberts and Srinivas~\cite{RS} proved the following:
Assume that the standard conjecture and Bloch-Beilinson conjecture are true.
Then $\phi$ is an isomorphism if the defining ideal of $R$ is generated by polynomials with coefficients in the algebraic closure of the prime field.
\end{itemize}
\end{itemize}
\end{rm}
\end{Example}

Consider the groups $\overline{G_0(R)} \subset \overline{G_0(R)}_{\Bbb Q} \subset  \overline{G_0(R)}_{\Bbb R}$.
We shall define some cones in $\overline{G_0(R)}_{\Bbb R}$.

\begin{Definition}
\begin{rm}
Let $C_{CM}(R)$ be the cone (in $\overline{G_0(R)}_{\Bbb R}$) spanned by all maximal Cohen-Macaulay $R$-modules.  
\[
C_{CM}(R) = \sum_{M: MCM}{\Bbb R}_{\ge 0}[M] \subset \overline{G_0(R)}_{\Bbb R} .
\]
We call it the {\em Cohen-Macaulay cone} of $R$.
Thinking a free basis of $\overline{G_0(R)}$ as an 
orthonormal basis of $\overline{G_0(R)}_{\Bbb R}$,
we think $\overline{G_0(R)}_{\Bbb R}$ as a metric space.
Let $C_{CM}(R)^-$ be the closure of $C_{CM}(R)$
with respect to this topology on $\overline{G_0(R)}_{\Bbb R}$.

We define the {\em strictly nef cone} by
\[
SN(R) = \{ \alpha \mid 
\mbox{ $\chi_{{\Bbb F}.}(\alpha) > 0$ for any ${\Bbb F}. \in C_d(R)$} \} .
\]
\end{rm}
\end{Definition}

By the depth sensitivity, $\chi_{{\Bbb F}.}([M]) = \ell_R(H_0({\Bbb F}.\otimes M))
> 0$ for any maximal Cohen-Macaulay module $M \ (\neq 0)$ and  ${\Bbb F}. \in C_d(R)$.   Therefore, 
\[
SN(R) \supset C_{CM}(R) - \{ 0 \} .
\]

\begin{Remark}\label{rem3}
\begin{rm}
Assume that $R$ is a Cohen-Macaulay local domain.
Let $M$ be a torsion $R$-module.
Taking sufficiently high syzygies of $M$, we know 
\[
\pm [M] + n [R] \in C_{CM}(R) \ \ \mbox{for $n \gg 0$}.
\]
Therefore, we have $\dim C_{CM}(R) = \rank \overline{G_{0}(R)}$ and
\[
C_{CM}(R)^{-} \supset C_{CM}(R) \supset Int(C_{CM}(R)^{-})
= Int(C_{CM}(R)) \ni [R] ,
\]
where $Int(~~)$ denotes the interior. 
\end{rm}
\end{Remark}

\begin{Example}
\begin{rm}
The following examples are given in \cite{DK2}.
Assume that $k$ is an algebraically closed field of characteristic zero.
\begin{itemize}
\item[(1)]
Put $R = k[x,y,z,w]_{(x,y,z,w)}/(xy-f_1f_2\cdots f_t)$.
Here, we assume that $f_1$, $f_2$, \ldots, $f_t$ are pairwise coprime
linear forms in $k[z,w]$ with $t \ge 2$.
In this case, we have $\rank \overline{G_{0}(R)} = t$.
We know (see \cite{DK2}) that the Cohen-Macaulay cone is minimally spanned by
the following $2^t - 2$ maximal Cohen-Macaulay modules of rank one:
\[
\{ (x, f_{i_1}f_{i_2}\cdots f_{i_s}) \mid
1 \le s < t, \ \ 1 \le i_1 < i_2 < \cdots < i_s \le t \}
\]

Here, remark that this ring is of finite representation type if and only if $t \le 3$.
\item[(2)]
The Cohen-Macaulay cone of 
$k[x_1,x_2, \ldots, x_6]_{(x_1,x_2, \ldots, x_6)}/(x_1x_2+x_3x_4+x_5x_6)$
is not spanned by maximal Cohen-Macaulay modules of rank one.
It is of finite representation type since it has a simple singularity.
\end{itemize}
\end{rm}
\end{Example}

\section{Fundamental class}

\begin{Definition}
\begin{rm}
Let $R$ be a $d$-dimensional Noetherian local domain.
We put 
\[
\mu_R = {\tau_R}^{-1}([\Spec R]) \in G_0(R)_{\Bbb Q} ,
\]
where $\tau_R : G_0(R)_{\Bbb Q} \stackrel{\sim}{\rightarrow} A_*(R)_{\Bbb Q}$ 
is the singular Riemann-Roch map, and
$[\Spec R]$ denotes the cycle in $A_*(R)$ corresponding to the scheme
$\Spec R$ itself.
\[
\begin{array}{ccc}
G_0(R)_{\Bbb Q} & \longrightarrow & \overline{G_0(R)}_{\Bbb Q} \\
\mu_R & \mapsto & \overline{\mu_R}
\end{array}
\]
We call the image of $\mu_R$ in $\overline{G_0(R)}_{\Bbb Q}$ 
the {\em fundamental class} of $R$, and denote it by $\overline{\mu_R}$.
\end{rm}
\end{Definition}

Remark that $\overline{\mu_R} \neq 0$ since $\rank_R \mu_R = 1$.

Put $R = T/I$, where $T$ is a regular local ring.
The map $\tau_R$ is defined using not only $R$ but also $T$.
Therefore,  $\mu_R$
may depend on the choice of $T$.\footnote{There is no example that the map $\tau_R$
actually depend on the choice of $T$.
For some excellent rings, it had been proved that $\tau_R$ is independent of
the choice of $T$ (Proposition~1.2 in \cite{K16}).}
However, we can prove that $\overline{\mu_R}$ is independent of $T$ (Theorem~5.1 in \cite{K23}). 

We shall explain the reason why we call $\overline{\mu_R}$ the fundamental class of $R$.

\begin{Remark}\label{Remark8}
\begin{rm}
\begin{itemize}
\item[(1)]
If $X$ ($= \Spec R$) is a $d$-dimensional affine variety over ${\Bbb C}$,
we have the cycle map $cl$
such that $cl([\Spec R])$ coincides with the fundamental class $\mu_X$ in 
$H_{2d}(X, {\Bbb Q})$ in the  usual sense,
where $H_*(X, {\Bbb Q})$ is the Borel-Moore homology.
Here $\mu_X$ is the generator of $H_{2d}(X, {\Bbb Q}) \simeq {\Bbb Z}$. 
\[
\begin{array}{ccccc}
G_{0}(R)_{\Bbb Q} & \stackrel{\tau_{R}}{\longrightarrow} & A_{*}(R)_{\Bbb Q} &
\stackrel{cl}{\longrightarrow} & H_{*}(X, {\Bbb Q}) \\
\mu_R & \mapsto & [\Spec R] & \mapsto & \mu_{X}
\end{array}
\]
The map $cl$ induces the map $\overline{A_{d}(R)}_{\Bbb Q} \longrightarrow H_{2d}(X, {\Bbb Q})$
such that the fundamental class $\mu_X$ is the image of $\overline{\tau_R}(\overline{\mu_R})$. 
Hence, we call $\overline{\mu_R}$ the fundamental class of $R$.
\item[(2)]
Let $R$ have a subring $S$ such that $S$ is a regular local ring and
$R$ is a localization of a finite extension of $S$.
Let $L$ be a finite-dimensional normal extension of $Q(S)$ containing $Q(R)$.
Let $B$ be the integral closure of $R$ in $L$.
Then, we have
\[
\mbox{$\mu_R = \frac{1}{\rank_RB}[B]$ in $G_0(R)_{\Bbb Q}$. }
\]
In particular, $\overline{\mu_R} = \frac{[B]}{\rank_RB}$ in $\overline{G_0(R)}_{\Bbb Q}$ (see the proof of Theorem~1.1 in \cite{K15}).
\item[(3)]
Assume that $R$ is of characteristic $p > 0$ and F-finite.
Assume that the residue class field is algebraically closed.
By the singular Riemann-Roch theorem, we have
\[
\mbox{${\displaystyle \overline{\mu_R} = \lim_{e \to \infty} \frac{[{}^eR]}{p^{de}}}$ in 
$\overline{G_0(R)}_{\Bbb R}$,}
\]
where ${}^eR$ is the $e$-th Frobenius direct image
(see Definition~\ref{directIm}, \ref{Ffin} below).
It immediately follows from the equations~(\ref{4}) and (\ref{last}) below.
\end{itemize}
\end{rm}
\end{Remark}

\begin{Example}\label{ex7}
\begin{rm}
\begin{itemize}
\item[(1)]
If $R$ is a complete intersection, then $\mu_R$ is equal to $[R]$ in 
 ${G_0(R)}_{\Bbb Q}$, therefore $\overline{\mu_R} = [R]$ in $\overline{G_0(R)}_{\Bbb Q}$.
There exists a Gorenstein ring such that $\overline{\mu_R} \neq [R]$.
However there exist many examples of rings satisfying $\overline{\mu_R} = [R]$
(\cite{K16}).
Roberts (\cite{R2}, \cite{R1}) proved the vanishing property of intersection multiplicities
for rings satisfying $\overline{\mu_R} = [R]$.
\item[(2)]
Let $R$ be a normal domain.
Then, we have
\[
\begin{array}{ccl}
G_{0}(R)_{\Bbb Q} & \stackrel{\tau_{R}}{\longrightarrow} & A_{*}(R)_{\Bbb Q} 
= A_{d}(R)_{\Bbb Q} \oplus A_{d-1}(R)_{\Bbb Q} \oplus \cdots \\
\mbox{$[R]$} & \mapsto & [\Spec R] -\frac{K_{R}}{2} + \cdots \\
\mbox{$[\omega_{R}]$} & \mapsto & [\Spec R] + \frac{K_{R}}{2} + \cdots ,
\end{array}
\]
where $K_R$ is the Weil divisor corresponding to the canonical module $\omega_R$.
If $\tau_{R}^{-1}(K_{R}) \neq 0$ in $\overline{G_{0}(R)}_{\Bbb Q}$,
then $[R] \neq \overline{\mu_R}$. 
Although the equality
\[
\overline{\mu_R} = \frac{1}{2}([R] + [\omega_{R}])
\]
is sometimes satisfied,
it is not true in general.
\item[(3)]
Let $R = k[x_{ij}]/I_{2}(x_{{ij}})$, where $(x_{ij})$ is the 
generic $(m+1) \times (n+1)$-matrix, and $k$ is a field.
Suppose $0 < m \le n$.
Then, we have
\[
\begin{array}{rcl}
G_{0}(R)_{\Bbb Q} \simeq \overline{G_{0}(R)}_{\Bbb Q} & \simeq &
A_*(R)_{\Bbb Q} \simeq {\Bbb Q}[a]/(a^{m+1}) \\
\mbox{$[R]$} & \mapsto & \left( \frac{a}{1-e^{-a}} \right)^{m}
\left( \frac{-a}{1-e^{a}} \right)^{n} \\
& & = 1 + \frac{1}{2}(m-n)a + \frac{1}{24}( \cdots )a^{2} + \cdots \\
\mbox{$[\omega_R]$} & \mapsto & \left( \frac{-a}{1-e^{a}} \right)^{m}
\left( \frac{a}{1-e^{-a}} \right)^{n} \\
\overline{\mu_R} & \mapsto & 1 \\
\tau_R^{-1}(K_{R}) & \mapsto & (n-m)a 
\end{array}
\]
\item[(4)]
By Remark~2.9 in \cite{CK}, if $\overline{\mu_R} \in C_{CM}(R)$,
then there exists a maximal Cohen-Macaulay $R$-module $N$ such that
$[N] = \rank_RN \cdot \overline{\mu_R}$ in $\overline{G_0(R)}_{\Bbb Q}$.
\end{itemize}
\end{rm}
\end{Example}

Here, we shall explain the connection between the fundamental class $\overline{\mu_R}$
and the homological conjectures.

\begin{Fact}\label{fact8}
\begin{rm}
\begin{itemize}
\item[(1)]
The small Mac conjecture is true if and only if  $\overline{\mu_R} \in C_{CM}(R)$ for any complete local domain $R$ (Theorem~1.3 in \cite{K15}).
We give an outline of the proof here.

``If" part is trivial.
We shall show ``only if" part.
Suppose that $S$ is a regular local ring such that
$R$ is a finite extension over $S$. 
Let $L$ be a finite-dimensional normal extension of $Q(S)$ containing $Q(R)$.
Let $B$ be the integral closure of $R$ in $L$.
Then, $B$ is finite over $R$, and $B$ is a complete local domain.
Here, {\em assume that there exists an maximal Cohen-Macaulay $B$-module $M$}.
Put ${\rm Aut}_{Q(S)}(L) = \{ g_1, \ldots, g_t \}$ and $N = \oplus_i ({}_{g_i}M)$,
where ${}_{g_i}M$ denotes $M$ with $R$-module structure given by 
$a \times m = g_i(a)m$.
Then $N$ is a maximal Cohen-Macaulay $R$-module such that 
$[N] = \rank_RN \cdot \mu_R$ in $G_0(R)_{\Bbb Q}$.
Therefore, $\overline{\mu_R} = \frac{[N]}{\rank_RN} \in C_{CM}(R)$.

{\em Even if $R$ is an equi-characteristic Gorenstein ring, it is not known whether
$\overline{\mu_R}$ is in $C_{CM}(R)$ or not}.
If $R$ is a complete intersection, then $\overline{\mu_R} = [R] \in C_{CM}(R)$
as in (1) in Example~\ref{ex7}.
\item[(2)]
If $\overline{\mu_R} = [R]$ in $\overline{G_0(R)}_{\Bbb Q}$, then 
the vanishing property of intersection multiplicities holds (Roberts~\cite{R2}, 
\cite{R1}). 

\item[(3)]
Roberts~\cite{Ri} proved $\overline{\mu_R} \in SN(R)$ if $ch(R) = p > 0$.
Using it, he proved the new intersection theorem 
in the mixed characteristic case.

\item[(4)]
If $R$ contains a field, then  $\overline{\mu_R} \in SN(R)$ 
(Kurano-Roberts~\cite{K13}).
{\em Even if $R$ is a Gorenstein ring (of mixed characteristic),
we do not know whether $\overline{\mu_R} \in SN(R)$ or not.}

\item[(5)]
If $\overline{\mu_R} \in SN(R)$ for any $R$,
then Serre's positivity conjecture is true in the case where
one of two modules is (not necessary maximal) Cohen-Macaulay.

It is well-known that Serre's positivity conjecture follows from the small Mac conjecture. 
\end{itemize}
\end{rm}
\end{Fact}

\begin{Remark}\label{rem9}
\begin{rm}
\begin{itemize}
\item[(1)]
If $R$ is Cohen-Macaulay of characteristic $p > 0$, then
${}^eR$ is a maximal Cohen-Macaulay module.
Since $\overline{\mu_R}$ is the limit of $[{}^eR]/p^{de}$ in $\overline{G_0(R)}_{\Bbb R}$ as in Remark~\ref{Remark8} (3),
$\overline{\mu_R}$ is contained in $C_{CM}(R)^-$.
If we know that $C_{CM}(R)$ is a closed set of $\overline{G_0(R)}_{\Bbb R}$,
we have  $\overline{\mu_R} \in C_{CM}(R)^- = C_{CM}(R)$.
If the cone $C_{CM}(R)$ is finitely generated, then it is a closed subset.
We do not know any example that the cone $C_{CM}(R)$ is not finitely generated.

{\em In the case where $R$ is not of characteristic $p > 0$,
we do not know whether $\overline{\mu_R}$ is contained in $C_{CM}(R)^-$
even if $R$ is a Gorenstein ring.}
\item[(2)]
As we have already seen in Remark~\ref{rem3}, if $R$ is Cohen-Macaulay, 
then $[R] \in Int(C_{CM}(R)) \subset C_{CM}(R)$.

There is an example of non-Cohen-Macaulay ring $R$ containing a field such that
$[R] \not\in SN(R)$.\footnote{It was conjectured above 50 years ago
that $[R]$ was in $SN(R)$ for any local ring $R$.
Essentially, the famous counter example 
due to Dutta-Hochster-MacLaughlin~\cite{DHM} gives an example
 $[R] \not\in SN(R)$.}
On the other hand, it is expected that $\overline{\mu_R} \in SN(R)$ for any $R$
(Fact~\ref{fact8} (4)).
Therefore, for a non-Cohen-Macaulay local ring $R$,
$\overline{\mu_R}$ behaves better than $[R]$ in a sense.
\end{itemize}
\end{rm}
\end{Remark}

\section{Main theorem}

In Fact~\ref{fact8}, we saw
that the fundamental class $\overline{\mu_R}$ is deeply related to the homological conjectures.
We propose the following question.

\begin{Question}\label{Q12}
\begin{rm}
Assume that $R$ is a ``good" Cohen-Macaulay local domain (for example, 
equi-characteristic, Gorenstein, etc).
Is $\overline{\mu_R}$ in $C_{CM}(R)$?
\end{rm}
\end{Question}

If $R$ is a Cohen-Macaulay local domain such that the rank of $\overline{G_{0}(R)}$ is one, 
then  $[R]= \overline{\mu_R} \in C_{CM}(R)$,
therefore Question~\ref{Q12} is true in this case.
There are a lot of such examples (for instance, invariant subrings with respect to finite group actions, etc.).

\begin{Definition}\label{directIm}
\begin{rm}
Let $p$ be a prime number and $R$ be a Noetherian ring of characteristic $p$. 
Let $e>0$ be an integer and 
\[
F^e:R \longrightarrow R
\]
be the $e$-th Frobenius map.
We denote by $^eR$ the $R$-module $R$ with $R$-module
structure given by $r \times x = F^e(r)x$. 
It is called the $e$-th {\em Frobenius direct image}.
\end{rm}
\end{Definition}

\begin{Definition}\label{Ffin}
\begin{rm}
Let $p$ be a prime number and $R$ be a Noetherian ring of characteristic $p$. 
We say that $R$ is {\em F-finite} if the Frobenius map $F: R \longrightarrow R$ is finite.
\end{rm}
\end{Definition}

\begin{Remark}\label{mitomeru}
\begin{rm}
Let $R$ be a $d$-dimensional F-finite Noetherian local ring. 
We have the following commutative diagram~(\ref{kakan}) where the horizontal map $\tau _R$ is the singular Riemann-Roch map and 
the vertical  maps are induced by $F^e$:
\begin{equation}\label{kakan}
 \begin{CD}
 G_0(R)_{\Bbb Q}          @>{\tau _R}>>          A_*(R)_{\Bbb Q}\\
 @V{F^e_*}VV                                              @VV{F^e_*}V\\
 G_0(R)_{\Bbb Q}         @>{\tau _R}>>           A_*(R)_{\Bbb Q}
 \end{CD}
\end{equation}\\
By diagram~(\ref{kakan}), we have
\begin{equation}\label{a}
\tau_R([^eR]) = F^e_*\bigl(\tau_R([R]) \bigr) .
\end{equation}
We set 
\[
\tau_R([R]) = \tau_R([R])_d + \tau_R([R])_{d-1} + \cdots + \tau_R([R])_{0}
\]
 where $\tau_R([R])_i \in A_i(R)_{\Bbb Q}$ for $i=0,\dots , d$. 
Then, by the top term property \cite{F}, we know
\begin{equation}\label{toptermproperty}
\tau_R([R])_d = [\Spec R] \in A_*(R)_{\Bbb Q} .
\end{equation}

Assume that $(R,{\mathfrak m})$ is a 
$d$-dimensional $\rm F$-$\rm finite$ Noetherian local domain with residue class field $R/{\mathfrak m}$ algebraically closed. 
For $\alpha \in A_i(R)_{\Bbb Q}$ we have 
\begin{equation}\label{egen}
F_*(\alpha)=p^i\alpha
\end{equation}
by Lemma~\ref{rank} below and the definition of $F_*$~\cite{F}.
Therefore
\begin{equation}\label{b}
F^e_*\bigl(\tau_R([R]) \bigr) = p^{de}[\Spec R] + \sum_{0 \le i \le d-1}p^{ie} \tau_R([R])_{i}
.
\end{equation}
Hence, by the equations $\eqref{a}$, $\eqref{b}$, we have
\[
\tau_R([{}^eR])_i = p^{ie} \tau_R([R])_i .
\] 
Therefore, 
\begin{equation}\label{4}
[ {}^eR]
=
p^{de} \tau^{-1}_R \bigl( [\Spec R] \bigr) + \sum_{0 \le i \le d-1}p^{ie} \tau^{-1}_R(\tau_R([R])_{i})  
\end{equation}
in $G_0(R)_{\Bbb Q}$. 
\end{rm}
\end{Remark}

The following lemma is well-known. We omit a proof. 

\begin{Lemma}\label{rank}
Assume that $R$ is an F-finite Noetherian local domain 
of characteristic $p$ with residue class field algebraically closed.
Then, for any $e>0$, we have
\[
\rank_R {}^eR  = p^{(\dim R)e} .
\]
\end{Lemma}

\begin{Definition}\label{FFRT}
\begin{rm}
Let $R$ be a Cohen-Macaulay ring of characteristic $p>0$. 
We say that $R$ is {\em FFRT} ({\em of  finite F-representation type}) if 
there exist finitely many indecomposable maximal Cohen-Macaulay $R$-modules $M_1 , \dots , M_s$ such that  
there exist nonnegative integers $a_{e1} , \dots , a_{es}$ with 
\[
{}^eR \simeq M^{a_{e1}}_1 \oplus \cdots \oplus M^{a_{es}}_s
\]
for each $e>0$.　
\end{rm}
\end{Definition}

\begin{Definition}\label{tightclosure}
\begin{rm}
Let $p$ be a prime number and $R$ be a Noetherian ring of characteristic $p$. 
Let $R^{\circ}$ be the set of elements of $R$ that are not contained in any minimal prime ideals of $R$. 
Let $I$ be an ideal of $R$.
Given a natural number $e$, set $q=p^e$. 
The ideal generated by the $q$-th powers of elements of $I$ 
is called the $q$-th Frobenius power of $I$, denoted by $I^{[q]}$. 
We define the {\em tight closure} $I^*$ of $I$ as follows:
\[
I^*=\{ x \in R \mid \mbox{there exists $c\in R^{\circ}$ such that $cx^q \in I^{[q]}$ for $q\gg 0$} \} .
\]
We say that $I$ is {\em tightly closed} if $I = I^*$.
\end{rm}
\end{Definition}

\begin{Definition}\label{F-rational}
\begin{rm}
Let $R$ be a Noetherian local ring of characteristic $p>0$. 
We say that $R$ is {\em F-rational} if every parameter ideal is tightly closed.
\end{rm}
\end{Definition}

\vspace{2mm}

Now, we start to prove Theorem~\ref{th11} (1).
Since $R$ is FFRT, 
there exist finitely many indecomposable maximal Cohen-Macaulay $R$-modules $M_1 , \dots , M_s$ such that  
there exist nonnegative integers $a_{e1} , \dots , a_{es}$ with 
\begin{equation}\label{bunkai}
{}^eR \simeq M^{a_{e1}}_1 \oplus \cdots \oplus M^{a_{es}}_s
\end{equation}
for each $e>0$. 
Let $U$ be the ${\Bbb Q}$-vector subspace of $G_0(R)_{\Bbb Q}$ spanned 
by 
\[
\{ [M_1], \dots , [M_s] \} \cup \{ \tau ^{-1}_R \bigl( \tau _R([R])_j  \bigr) \mid 0 \le j \le d \} .
\]
Here, recall that $\mu_R = \tau_R^{-1}(\tau_R([R])_d) \in U$
by the top term property (\ref{toptermproperty}).
Although we can show that $U$ is spanned by $\{ [M_1], \dots , [M_s] \}$, 
we do not need it in this proof.
Thinking a basis of $U$ as an 
orthonormal basis of $U_{\Bbb R}$,
we think $U_{\Bbb R}$ as a metric space.
Set $C= \displaystyle\sum^{s}_{i=1}{\Bbb R}_{\ge 0}[M_i] \subset U_{\Bbb R}$. 
Then $C$ is a closed subset of $U_{\Bbb R}$.
We shall show $\mu_R \in C$. 

Since the residue field is algebraically closed, 
$\rank_R {}^eR=p^{de}$ for any $e>0$ by Lemma~\ref{rank}. 
Since 
\[
[{}^e R] = a_{e1}[M_1] + \cdots + a_{es}[M_s] 
\]
by (\ref{bunkai}),
we have 
\[
\displaystyle\frac{1}{p^{de}}[{}^eR] \in C 
\] 
for any $e > 0$.
By the equation~(\ref{4}), 
\begin{equation}\label{last}
\displaystyle\frac{1}{p^{de}} [{}^eR] = \displaystyle\sum_{0 \le i \le d} \frac{1}{p^{ie}}\tau^{-1}_R \bigl( \tau_R([R])_{d-i} \bigr)
.
\end{equation}
By the definition of $U$, every term of the right-hand side is in $U_{\Bbb R}$. 
Hence we have 
\[
\displaystyle\lim_{e \to \infty} \frac{1}{p^{de}}[{}^eR]  = \tau^{-1}_R \bigl( \tau_R([R])_d \bigr) = \tau^{-1}_R \bigl( [\Spec R] \bigr) 
={\mu_R} \ \ \mbox{in $U_{\Bbb R}$}.
\] 
Since $C$ is a closed set of $U_{\Bbb R}$, we have
$\mu_R \in C$.
By the same argument as in Example~\ref{ex7} (4),
there exist a natural number $n$ and a 
maximal Cohen-Macaulay $R$-module $N$ such that
$n \mu_R = [N]$ in $G_0(R)_{\Bbb Q}$.

\vspace{2mm}

Next, we start to prove Theorem~\ref{th11} (2).

First, we shall prove that $[\omega_R] \in Int(C_{CM}(R))$
if $R$ is Cohen-Macaulay.
We have a homomorphism $\xi : {G_0(R)}_{\Bbb R} \rightarrow
{G_0(R)}_{\Bbb R}$ given by 
$\xi([M]) = \sum_i(-1)^i[{\rm Ext}_R^i(M,\omega_R)]$.
For a maximal Cohen-Macaulay module $M$,
${\rm Ext}^i_R(M,\omega_R) = 0$ for $i > 0$ and
${\rm Hom}_R({\rm Hom}_R(M,\omega_R),\omega_R) \simeq M$.
Therefore, $\xi^2$ is equal to the identity, and
$\xi$ is an isomorphism.
By the definition of $\tau_R$,
we have a commutative diagram\footnote{
Put $R = T/I$, where $T$ is a regular local ring.
Then, $\xi([M]) = (-1)^{{\rm ht}(I)} \sum_i(-1)^i[{\rm Ext}_T^i(M,T)]$. 
Let ${\Bbb F}.$ be a $T$-free resolution of $M$.
Then, by the definition of $\tau_R$, we have $\tau_R([M]) =
{\rm ch}({\Bbb F}.) \cap [\Spec T]$, where ${\rm ch}({\Bbb F}.)$
is the localized Chern character of ${\Bbb F}.$ (\S 18 in \cite{F}).
By the local Riemann-Roch formula (Example~18.3.12 in \cite{F}),
$\tau_{R}(\xi([M])) = {\rm ch}({\Bbb F}.^*[{\rm ht}(I)]) \cap [\Spec T]$.
By Example~18.1.2, we obtain the equality~(\ref{altsum}). }
\[
\begin{array}{ccc}
 {G_0(R)}_{\Bbb R} & \stackrel{\tau_R\otimes 1}{\longrightarrow} &
A_*(R)_{\Bbb R} \\
{\scriptstyle \xi}\downarrow \hphantom{\scriptstyle \xi} & & 
{\scriptstyle \phi}\downarrow \hphantom{\scriptstyle \phi} \\
 {G_0(R)}_{\Bbb R} & \stackrel{\tau_R\otimes 1}{\longrightarrow} &
A_*(R)_{\Bbb R}
\end{array}
\]
where $\phi : A_*(R)_{\Bbb R} \rightarrow A_*(R)_{\Bbb R}$ is the map
given by
\begin{equation}\label{altsum}
\phi(q_d + q_{d-1} + \cdots + q_{i} + \cdots + q_0)
= q_d - q_{d-1} + \cdots +(-1)^{d-i} q_{i} + \cdots + (-1)^dq_0
\end{equation}
for $q_i \in A_i(R)_{\Bbb R}$.
Since the numerical equivalence is graded in $A_*(R)_{\Bbb Q}$ 
as in Proposition~2.4 in \cite{K23},
$\phi$ preserves the numerical equivalence.
Therefore we have the induced map
\[
\overline{\xi} : \overline{G_0(R)}_{\Bbb R} \rightarrow
\overline{G_0(R)}_{\Bbb R} .
\]
Remark that $\overline{\xi}$ is an isomorphism of $\R$-vector spaces since 
$\overline{\xi}^2$
is the identity. 
The map $\overline{\xi}$ satisfies $\overline{\xi}([R]) = [\omega_R]$ and $\overline{\xi}(C_{CM}(R)) = C_{CM}(R)$.
Since $[R] \in Int(C_{CM}(R))$ by Remark~\ref{rem3},  we obtain $[\omega_R] \in Int(C_{CM}(R))$.

Assume that $M$ is a maximal Cohen-Macaulay module.
For $e > 0$, consider the following exact sequence
\[
0 \longrightarrow L_e \longrightarrow F^e_{*}(M)
 \longrightarrow M^{\oplus b_e} \longrightarrow 0
\]
where $F^e_{*}(M)$ is the $e$-th Frobenius direct image of $M$.
Take $b_e$ as large as possible.
Recall that $L_e$ is a maximal Cohen-Macaulay module.
Put $r = \rank_R M$.

Here we define the dual F-signature following Sannai~\cite{Sa} as follows: 
\[
s(M) = \limsup_{e \to \infty} \frac{b_e}{r p^{de}} 
\]
Then, taking a subsequence of $\{ \frac{b_e}{rp^{de}} \}_e$,
we may assume that $s(M) = \lim_{e \to \infty} \frac{b_e}{r p^{de}}$.

On the other hand, consider
\[
\tau_R([M]) = \tau_R([M])_d + \tau_R([M])_{d-1} + \cdots + \tau_R([M])_{0}.
\]
Here, we have $ \tau_R([M])_d = r[\Spec R]$
since $[M] - r[R]$ is a sum of cycles of torsion modules.
By (\ref{kakan}) and (\ref{egen}), 
\begin{eqnarray*}
\tau_R([F^e_*(M)]) & = & F^e_*( \tau_R([M])_d + \tau_R([M])_{d-1} + \cdots + \tau_R([M])_{0}) \\
& = & p^{de}\tau_R([M])_d + p^{(d-1)e}\tau_R([M])_{d-1} + \cdots + \tau_R([M])_{0} .
\end{eqnarray*}
Then, we have
\[
\overline{\tau_R}(\lim_{e \to \infty} \frac{[F^e_*(M)]}{rp^{de}}) =\frac{ \tau_R([M])_d}{r} = [\Spec R] \ \ \mbox{in $\overline{A_*(R)}_{\Bbb R}$.}
\]
Thus, 
\[
\lim_{e \to \infty} \frac{[F^e_*(M)]}{rp^{de}} = \overline{\mu_R}
 \ \ \mbox{in $\overline{G_0(R)}_{\Bbb R}$.}
\]
Then, $\frac{[L_e]}{rp^{de}}$ converges to some element in 
$\overline{G_0(R)}_{\Bbb R}$, say $\alpha(M)$.
\[
\begin{array}{cccccl}
\frac{[F^e_{*}(M)]}{rp^{de}} & = & \frac{b_e[M]}{rp^{de}} & + &
\frac{[L_e]}{rp^{de}} & \in  \overline{G_0(R)}_{\Bbb R} \\
\downarrow & & \downarrow & & \downarrow &  {\scriptstyle (e \to \infty)} \\
\overline{\mu_R} & = & s(M)[M] & + & \alpha(M) &  
\end{array}
\]
Since $L_e$ is  a maximal Cohen-Macaulay module,
we know $\alpha(M) \in C_{CM}(R)^-$.

Here set $M = \omega_R$.
Then
\begin{equation}\label{1}
\overline{\mu_R}  =  s(\omega_R)[\omega_R]  +  \alpha(\omega_R) \in 
\overline{G_0(R)}_{\Bbb R} ,
\end{equation}
where
\begin{equation}\label{2}
\alpha(\omega_R) \in C_{CM}(R)^- 
\end{equation}
and
\begin{equation}\label{3}
[\omega_R] \in Int(C_{CM}(R)) = Int(C_{CM}(R)^{-}) .
\end{equation}

The most important point in this proof is the fact that
\[
\mbox{$R$ is F-rational if and only if $s(\omega_R) > 0$}
\]
due to Sannai~\cite{Sa}.

Therefore, if $R$ is F-rational, then $\overline{\mu_R} \in Int(C_{CM}(R))$
by (\ref{1}),  (\ref{2}), (\ref{3}) and Remark~\ref{rem3}.
\qed

\begin{Remark}
\begin{rm}
If $R$ is a toric ring (a normal semi-group ring over a field $k$),
then we can prove $\overline{\mu_R} \in C_{CM}(R)$ as in the case FFRT
without assuming that $ch(k)$ is positive.
\end{rm}
\end{Remark}

\begin{Problem}
\begin{rm}
\begin{itemize}
\item[(1)]
As in the above proof, if there exists a maximal Cohen-Macaulay module
 in $Int(C_{CM}(R))$ such that
its generalized F-signature or its dual F-signature is positive,
then $\overline{\mu_R}$ is in  $Int(C_{CM}(R))$.

Without assuming that $R$ is F-rational,
do there exist such a maximal Cohen-Macaulay module?
\item[(2)]
How do we make mod $p$ reduction? \ \ (for example, the case of rational singularity)
\item[(3)]
If $R$ is Cohen-Macaulay,
is $\overline{\mu_R}$ in $C_{CM}(R)^{-}$?
If $R$ is a Cohen-Macaulay ring containing a field of positive characteristic,
then $\overline{\mu_R}$ in $C_{CM}(R)^{-}$ as in (1) in Remark~\ref{rem9}.
\item[(4)]
If $R$ is of finite representation type,
is $\overline{\mu_R}$ in $C_{CM}(R)$?
\item[(5)]
Find more examples of $C_{CM}(R)$ and $SN(R)$.
\end{itemize}
\end{rm}
\end{Problem}

In order to prove the following corollary,
it is enough to construct a $d$-dimensional Cohen-Macaulay local domain $A$
satisfying the following two conditions (Lemma~3.1 in \cite{CK}):
\begin{itemize}
\item[(1)]
$\overline{A_i(A)} \neq 0$ for $d/2 < i \le d$, and
\item[(2)] 
$\overline{\mu_A}$ is contained in $Int(C_{CM}(A))$.
\end{itemize}
The ring $R$ in Corollary~\ref{cor14} is the idealization of $A$ and certain maximal
Cohen-Macaulay $A$-module $M$.
We can simplify the proof of Corollary~\ref{cor14} using Theorem~\ref{th11}.
We know that  $k[x_{ij}]_{(x_{ij})}/I_2(x_{ij})$
satisfies the conditions (1) and (2) above,
where $(x_{ij})$ is the generic $n \times n$ or $n \times (n+1)$ matrix, and $I_2(x_{ij})$
stands for the ideal generated by $2$-minors of $(x_{ij})$.
In fact, by Example~\ref{rei} (2) (b) and Example~\ref{ex7} (3), the condition~(1) is satisfied.
Since  $k[x_{ij}]_{(x_{ij})}/I_2(x_{ij})$ is F-rational,  the condition~(2) is satisfied by Theorem~\ref{th11} (2).

\begin{Corollary}[\cite{CK}]\label{cor14}
Let $d$ be a positive integer and $p$ a prime number.
Let $\epsilon_0$, $\epsilon_1$, \ldots, $\epsilon_d$ be integers such that
\[
\epsilon_i = 
\left\{
\begin{array}{ll}
1 & i = d , \\
\mbox{$-1$, $0$ or $1$} & d/2 < i < d , \\
0 & i \le d/2 .
\end{array}
\right.
\]

Then, there exists a $d$-dimensional Cohen-Macaulay local ring $R$ of
characteristic $p$, a maximal primary ideal $I$ of $R$ of finite projective dimension,
and positive rational numbers $\alpha$, $\beta_{d-1}$, $\beta_{d-2}$,\ldots, $\beta_{0}$ such that
\[
\ell_R(R/I^{[p^n]}) = \epsilon_d \alpha p^{dn} + \sum_{i = 0}^{d-1} \epsilon_i \beta_i p^{in}
\]
for any $n > 0$.
\end{Corollary}

\noindent
Department of Mathematics\\
School of Science and Technology\\
Meiji University\\
Higashimata 1-1-1, Tama-ku \\
Kawasaki 214-8571, Japan

\end{document}